\documentclass{article}
\usepackage{latexsym}
\newtheorem{lemma}{LEMMA}
\newtheorem{nim}{THEOREM}
\newtheorem{cor}{COROLLARY}
\begin{document}
\title{How Far Can Nim in Disguise be Stretched?}
\author{Uri Blass \\
Electrical Engineering\\
Tel Aviv University\\
Ramat Aviv 61391, Israel
\and Aviezri S.\ Fraenkel \\
Applied Mathematics and Computer Science\\ 
Weizmann Institute of Science\\ 
Rehovot 76100, Israel
\and Romina Guelman\\
Institute of Mathematics\\ 
Hebrew University of Jerusalem\\ 
Jerusalem 91904, Israel}
\maketitle
\begin{abstract}
A move in the game of nim consists of taking any positive number of 
tokens from a {\sl single\/} pile. Suppose we add the class of moves 
of taking a nonnegative number of tokens jointly from {\sl all\/} the 
piles. We give a complete answer to the question which moves in the 
class can be adjoined without changing the winning strategy of nim. The 
results apply to other combinatorial games with unbounded Sprague-Grundy 
function values. We formulate two weakened conditions of the notion of 
nim-sum 0 for proving the results.  
\end{abstract}
\section{Introduction}

A cardinal theme in the theory of combinatorial games is how 
to generate new games from a given game or from a restricted 
class of games. The most widely used method is that of producing a 
game which is the {\it sum\/} of given games, but there are several 
other, less well-known methods; see e.g., chapter~14 of [Con1976]. 

A typical game consists of a finite collection of piles of finitely 
many tokens, where the moves are to remove a positive number of 
tokens from any single pile, or a positive number from several piles, 
according to specified rules. Such games often have equivalent 
manifestations, say in the form of board games, but for concreteness 
we shall restrict attention to the former. 

A central role in such games is taken by the game of {\sl nim\/}, 
in which only removal from any single pile is permitted. Most of our 
discussions will be centered about nim, but actually our results hold 
for any game which has unbounded Sprague-Grundy function values. 
Basic facts on the theory of combinatorial games can be found, e.g., 
in [BCG1982], [Con1976], [Guy1991], [Now1996]. 

Recently we began investigating the generation of new games by adding 
to given games classes of new moves [FrL1991], [FrO1998]. 
For a brief expository description of this approach, see 
[Fra1996, \S 6]. To conduct this program in an efficient way, 
it is very useful to find first the precise class of moves that can 
be adjoined to nim {\sl without\/} changing its winning strategy. 
This then allows to adjoin moves for which we will know that they 
modify the strategy of nim. 

Fundamental to investigations in combinatorial game theory is the 
notion of nim-sum. Let $S=\{a_1,\dots, a_n\}$ be a multiset set of 
nonnegative integers which has some 1-bit in a least significant position 
$k$, so to the right of position $k$ all the $a_i$ have 0-bits only. 
(Note that $S$ is a multiset rather than a set; the $a_i$ are 
not necessarily distinct.) 
If $S$ has nim-sum $\sigma=0$, we also say that $S$ is {\it even\/}, 
since every column in the binary expansions of $a_1,\dots, a_n$ 
has then an even number of 1-bits. We define $S$ to be {\it baLanced\/}, 
if $\sigma$ has a 0 in position $k$, since then the {\it L}east significant 
binary position in which the $a_i$ have 1-bits has an even number of 
1-bits. For $S$ to be balanced, no parity requirements are imposed on 
any digital position 
to the left of $k$. Finally, we say that $S$ is {\it smooTh\/} if it 
is balanced and $\sigma$ has a 0 also in position $k+1$, since then the 
last {\it T}wo binary positions in which the $a_i$ have 1-bits have an even 
number of 1-bits. Note that every even multiset is smooth, and 
every smooth multiset is balanced. A balanced multiset is a weaker 
form of a smooth multiset, and a smooth multiset is a weaker form of 
an even multiset. 

Let $S=\{a_1,\dots, a_n\}$ $(n\ge 2)$ be a multiset of nonnegative integers 
with at least two distinct $a_i>0$. Let $\Gamma$ be a game consisting 
of $n$ piles of finitely many tokens where $a_1,\dots, a_n$ tokens 
can be removed from the $n$ piles (in addition to the option of removing 
any positive number of tokens from any {\sl single\/} pile, as mentioned 
above). The player making the last move wins, and the opponent loses. 

Let $t_1,\dots, t_n$ be nonnegative integers with $t_i\ge a_i$ for 
all $i$. For any integer $s$, write $s(b)$ for the binary representation 
of $s$, and \mbox{$\sum^{'}$} for nim-summation. We distinguish 
between three cases.
\begin{enumerate}
\item \mbox{$\sum_{i=1}^{'n} t_{i}(b)\ne\sum_{i=1}^{'n} (t_{i}(b)-a_{i}(b))$} 
for all $t_i\ge a_i$ $(i\in\{1,\dots, n\})$. Then the strategy of 
$\Gamma$ is that of nim. This is so, because the nim-sum of position 
$(t_1\dots, t_n)$, which is \mbox{$\sum_{i=1}^{'n} t_{i}(b)$}, is 
distinct from the nim-sum of its follower $(t_1-a_1,\dots, t_n-a_n)$, 
which is \mbox{$\sum_{i=1}^{'n} (t_i(b)-a_i(b))$}. Hence the nim-sum 
is the Sprague-Grundy function of $\Gamma$. In terms of the game-graph 
of $\Gamma$, the move options of removing $a_1,\dots, a_n$ are  
equivalent to new edges in this digraph between vertices of 
{\sl distinct\/} S-G function values. 
\item 
\mbox{$\sum_{i=1}^{'n} t_{i}(b)=\sum_{i=1}^{'n} (t_{i}(b)-a_{i}(b))=R$} 
for some $t_i\ge a_i$ $(i\in\{1,\dots, n\})$. In this case there is a 
``short-circuiting'' of the S-G function value $R$ of $\Gamma$. If $R$ 
is never 0, then the strategy of $\Gamma$ is still the same as that of 
nim, but if $\Gamma$ is a component in a sum with another game, say with 
S-G function value $R$, then the move option of removing $a_1,\dots, 
a_n$ {\sl does\/} change the strategy of this sum. 
\item \mbox{$\sum_{i=1}^{'n} t_{i}(b)=\sum_{i=1}^{'n} (t_{i}(b)-a_{i}(b))=0$} 
for some $t_i\ge a_i$ $(i\in\{1,\dots, n\})$. In this case a 0 of the 
S-G function is short-circuited, so the strategy of $\Gamma$ is necessarily
different from that of nim. 
\end{enumerate}

It is easy to see that case~3 holds if $S$ is an even multiset. In 
Theorem~1 we prove that case~2 holds if and only if $S$ is a balanced 
multiset. In Theorem~2 we give necessary and sufficient conditions 
for the stronger case~3 to hold. It turns out that the condition of 
a balanced multiset has to be strengthened only slightly for case~3 
to hold. 

The precise forms of Theorems~1 and 2 are formulated in \S 2. Proofs 
are given in \S 3. The proof of Theorem~1 is constructive; it provides 
an algorithm for producing the integers $t_1,\dots, t_n$ such that 
\mbox{$\sum_{i=1}^{'n} t_{i}(b)=\sum_{i=1}^{'n} (t_{i}(b)-a_{i}(b))$}.
Similarly for Theorem~2. 

\section{The Main Results}

It is useful to preface the following definition before stating our 
first result.

\textbf{Definition 1.}\quad Let $S=\{a_1,\dots,a_n\}$ be a multiset of 
nonnegative integers. Denote by $\sigma$ the nim-sum of the $a_i$. 
Let $k$ be the maximum integer such that $2^k|a_i$ for every 
$i\in\{1,\ldots,n\}$. If $\sigma^k=0$ (the bit in position $k$ of 
$\sigma$), then $S$ is {\it balanced}. Otherwise it is 
{\it imbalanced}. If $\sigma^k=\sigma^{k+1}=0$, then $S$ is {\it smooth}. 
If $\sigma=0$, then $S$ is {\it even}.\medskip

Note that position $k$ is the least significant position in 
which any of the $a_i$ has a 1-bit, so $a_i(b)^j=0$ for all $j<k$. 
For example, $\{2,3,4\}$ is imbalanced ($k=0$), $\{1,2,5\}$ is balanced 
($k=0$) but not smooth, $\{2,3,5\}$ is smooth but not even, and 
$\{1,2,3\}$ is even.\medskip

If $S=\{a_1,\dots,a_n\}$ is an even multiset, then case~3 holds, 
since it holds, in fact, for $t_i=a_i$. A special case is when all the 
$a_i$ are the same and $n$ is even, in which case even 
\mbox{$\sum_{i=1}^{'n} la=\sum_{i=1}^{'n} (l-1)a=0$} for every positive 
integer $l$. Since the notions of balanced and smooth multisets are weak 
forms of that of even multisets, we may expect a weaker result for the 
former. This is indeed the case; the interesting point is that the result 
is not all that weaker. 

\begin{nim}$\!${\bf .}
Let $S=\{a_1,\dots,a_n\}$ be a multiset of nonnegative integers, $n\ge 2$, 
with at least two $a_i>0$. Then there are integers 
$t_{1},\ldots,t_{n}$ with $t_{i}\geq a_{i}$ for all $i$, 
such that
\begin{equation}
\mbox{$\sum_{i=1}^{'n} t_{i}(b) = \sum_{i=1}^{'n} (t_{i}(b)-a_{i}(b))$}
\end{equation} 
if and only if $S$ is a balanced multiset. 
\end{nim}

The proof that if $S$ is imbalanced then there are no integers $t_i$ 
satisfying (1) was already given in [FrL1991], where the truth of the 
opposite direction was conjectured. Since the known direction is the 
easy one, and in order for this paper to be self-contained, we repeat 
the short proof below. 

Our second theorem gives necessary and sufficient conditions for the 
stronger result (case~3 above) to hold. It turns out that though $S$ 
even is certainly a sufficient condition, it is by no means necessary. 

\begin{nim}$\!${\bf .}
Let $S=\{a_1,\dots,a_n\}$ be as in Theorem~1. Then there are integers 
$t_{1},\ldots,t_{n}$ with $t_{i}\geq a_{i}$, $i\in\{1,\ldots,n\}$, such that
\begin{equation}
\mbox{$\sum_{i=1}^{'n} t_{i}(b) = \sum_{i=1}^{'n} (t_{i}(b)-a_{i}(b))=0$}
\end{equation}
if and only if either 
\begin{enumerate}
\item $n$ is odd and $S$ is balanced. 
\item $n$ is even, and: either $S$ is balanced and there is 
$i\in\{1,\dots, n\}$ such that $a_i(b)^k=0$ $($where $k$ is as in 
Definition~1$)$; or $S$ is smooth and $n\ge 4$; or $S$ is even. 
\end{enumerate}
\end{nim}

We then have, 
\begin{cor}$\!${\bf .}
For $n=2$, $(2)$ holds if and only if $S$ is even, if and only if 
$a_1=a_2$. 
\end{cor}

To summarize, adjoining the moves of removing $a_1,\dots, a_n$ from 
the piles results in a game with a strategy different from nim if and 
only if (2) is satisfied, which, for $n=2$, is equivalent to $a_1=a_2$. 
If only (1) is satisfied, then the resulting game has the same strategy 
as nim, but the strategy will be different if the game is a 
component in a sum of games. 

The new notions in this paper are those of balanced and smooth sets, 
which are weakened conditions of the notion of nim-sum 0.

\section{Proofs}

{\bf Notation}
\begin{enumerate}
\item For any real number $x$, denote by $\lfloor x\rfloor$ the largest 
integer $\le x$. 
\item For any positive integer $s$, denote by $s(b)=\sum_{j=0}^{m} s^j2^j$ 
the {\it binary representation\/} of $s$, where 
$m=\lfloor\log_2 s\rfloor$, and $s^j\in\{0,1\}$ for all $j$. 
\item Whenever we add nonnegative integers, say $a_1,\dots, a_n$, we 
put
\begin{displaymath} 
m=\max(\lfloor\log_2 a_1\rfloor, \dots, \lfloor\log_2 a_n\rfloor),
\end{displaymath} 
which is consistent with $m$ in 2. 
\item \mbox{$\sum^{'}$} and $\oplus$ denote nim-summation. 
\end{enumerate}
Note that for any positive integers $a$ and $d$, $a(b)+d(b)=(a+d)(b)$.\medskip

\textbf{Definition 2.}\quad In the (binary) addition $a(b)+d(b)$, there is a 
{\it carry integer\/} 
$c(b)$, where 
$c(b)^{j+1}$ is the carry-bit generated by $a(b)^j+d(b)^j+c(b)^j$, to 
be added to $a(b)^{j+1}+d(b)^{j+1}$, namely, $c(b)^{j+1}=1$ if 
$a(b)^j+d(b)^j+c(b)^j>1$, and $c(b)^{j+1}=0$ otherwise, where $c(b)^0=0$ 
and $j\in\{0,\ldots,m\}$; $m$ as in Notation~3.\medskip 

The addition rule, based on Definition~2, is summarized in 
Table~1.\medskip

\centerline{\textbf{Table 1}}\bigskip
\begin{center}
\begin{tabular}{ccc|cc}
$a(b)^j$&$c(b)^j$&$d(b)^j$&$\big(a(b)+d(b)\big)^j$&$c(b)^{j+1}$\\
\noalign{\vspace{5pt}}
\hline
\noalign{\vspace{5pt}}
0&0&0&0&0\\ 0&0&1&1&0\\ 0&1&0&1&0\\ 0&1&1&0&1\\ 1&0&0&1&0\\ 1&0&1&0&1\\
1&1&0&0&1\\ 1&1&1&1&1
\end{tabular}
\end{center}

\begin{lemma}$\!${\bf .}
Let $a$ and $d$ be two integers. Then, in the above notation, $a(b)+d(b)=
a(b)\oplus d(b)\oplus c(b)$, where $c(b)$ is the carry integer of 
$a(b)+d(b)$.
\end{lemma}
\textbf{Proof.}  
The sum $a(b)+d(b)$ is given in the 4-th column of Table~1. We see 
that it has a 1-bit precisely in those rows in which the first 3 
columns have an odd number of 1-bits, i.e., precisely in rows in which 
$a(b)\oplus d(b)\oplus c(b)=1$.\quad$\Box$\medskip

\textbf{Proof of Theorem~1.} Let $d_{i}=t_{i}-a_{i}$. Then (1) holds 
if and only if 
\begin{equation}
\mbox{$\sum_{i=1}^{'n} (a_{i}(b)+d_{i}(b))=\sum_{i=1}^{'n} d_{i}(b)$.}
\end{equation}
It thus suffices to examine under what conditions 
$d_{1},\ldots,d_{n}$ can be constructed such that (3) holds. 

By Lemma 1, for every $i\in\{1,\ldots,n\}$,
$a_{i}(b)+d_{i}(b)=a_{i}(b)\oplus d_{i}(b)\oplus c_{i}(b)$, where $c_{i}(b)$ 
is the carry integer of the sum of $a_{i}(b)$ and $d_{i}(b)$. Substituting 
into (3), we get $\sum_{i=1}^{'n}(a_{i}(b)\oplus d_{i}(b)\oplus c_{i}(b))=
\sum_{i=1}^{'n} d_{i}(b)$. Thus (3) holds if and only if 
\begin{equation}
\mbox{$\sum^{'n}_{i=1} (a_{i}(b)\oplus c_{i}(b)) = 0$}.
\end{equation}

In every position $<k$, $a_i(b)$ has no 1-bits for all $i$, hence 
in every position $\le k$, $c_i(b)$ has no 1-bits for all $i$, where 
$k$ is as in Definition~1. Thus if $S$ is imbalanced, then in position $k$ 
there is an odd number of 1-bits, so (4) cannot hold. Hence there are 
no integers $t_1,\dots, t_n$ satisfying (1). 

So from now on we can assume that $S$ is balanced. To construct 
$d_{1},\dots, d_n$ satisfying (3) we first construct 
$c_{1},\dots,c_{n}$ satisfying (4), in Algorithm NotNimdi below, and 
then show how to construct the $d_i$.

Given an integer $a(b)$, an integer $c(b)$ can be a carry integer of 
the sum of $a(b)$ with an unknown integer $d(b)$, if the following
{\it carry rules\/} are kept. These rules follow immediately from 
Definition~2. 
\begin{enumerate}
\item If $l$ is the rightmost 1-bit of $a(b)$, then for every 
$j<l$ we have $c(b)^{j+1}=0$. 

For $j\ge l$, we have: 
\item If $a(b)^{j}=c(b)^{j}=0$, then $c(b)^{j+1}=0$.
\item If $a(b)^{j}=c(b)^{j}=1$, then $c(b)^{j+1}=1$.
\item If $a(b)^{j}+c(b)^{j}=1$, then $c(b)^{j+1}\in \{0,1\}$.

Indeed, in case~4 we clearly have $c(b)^{j+1}=d(b)^j$.
\end{enumerate}

Let now $m=\max(\lfloor\log_2 a_1\rfloor ,\dots,\lfloor\log_2 a_n\rfloor)$. 
Note that even if every $d_{i}(b)$ has its leftmost 1-bit in a position 
$\le m$, i.e., $d_i<2^{m+1}$ for all $i\in\{1,\dots,n\}$, 
any carry integer $c_{i}(b)$ may still have a 1-bit in position $m+1$. 

Consider the $2n\times (m+2)$ matrix $M$ consisting of 
$a_{1}(b),\ldots,a_{n}(b)$ 
with a blank line after each $a_{i}(b)$, where the carry $c_{i}(b)$ will 
be constructed in Algorithm NotNimdi1 below. Because of the anomaly, in 
English, of writing from left to right, yet writing numbers with their 
significance increasing from right to left, we will number the columns 
of $M$, contrary to the common convention, from right $(0)$ to left 
$(m+1)$. Also the carry-bits will be constructed from position (column) 
$0$ to $m+1$.\medskip 

The following are the guidelines the algorithm will follow.

\begin{enumerate}
\item[A.] In every column of $M$, the number of 1-bits is even, which is 
necessary to satisfy (4). 
\item[B.] Every $c_{i}$ is constructed to be consistent with the above 
carry rules. 
\item[C.] For every $j\in\{k,\ldots,m+1\}$ there are $h,l\in\{1,\ldots,n\}$, 
$h\ne l$, such that $a_{h}(b)^{j}+c_{h}(b)^{j}=a_{l}(b)^{j}+c_{l}(b)^{j}=1$, 
where $k$ is as in Definition~1. 
\end{enumerate}

Property C is needed to ensure that A and B can be realized in every 
column of $M$. Indeed, suppose the $(j-1)$-th column of $M$ is the 
0-vector, and the $j$-th column contains a single 1-bit. Then there is 
no way of mending the $j$-th column to have an even number of 1-bits, 
as needed for consistency with the carry rules.\medskip

Note that if $S$ is balanced, then column $k$ of $M$ contains an 
{\sl even positive\/} number of 1-bits, and all columns to the right 
of $k$ are the 0-vector, provided that $c_i(b)^j=0$ for all 
$j\in\{0,\dots,k\}$, $i\in\{1,\dots,n\}$. This indeed holds 
by carry rule~1.\bigskip 

Suppose that the $(j-1)$-th position was constructed satisfying the 
above guidelines, and now the $j$-th position must be constructed. First, 
to satisfy the carry rules, if $a_{i}(b)^{j-1}=c_{i}(b)^{j-1}=1$, then 
we must put $c_{i}(b)^j=1$. Secondly, if the number of 1-bits in the 
$j$-th column is even but C is violated, then it has to be restored 
so as to leave the number of 1-bits even. Finally, if the number of 
1-bits in the $j$-th position is odd, then the algorithm must change 
it to even such that C is also satisfied. These requirements are 
reflected in Algorithm NotNimdi1 below. The word ``Nimdi'' was coined 
in [FrL1991]; it stands for {\it NIM\/} in {\it DI}$\,$sguise. Since 
in the present case we have balanced multisets, for which the moves 
may result in a non-nim strategy, the designation NotNimdi for the 
algorithm seemed appropriate.\bigskip

\begin{center}{\bf Algorithm\quad NotNimdi1}
\end{center}
\begin{enumerate}
\item For $j\leq k$, put $c_{i}(b)^j=0$ for all $i$.

\item For $j$ from $k+1$ to $m+1$ do:
\begin{enumerate}
\item For every $i\in\{1,\ldots,n\}$ for which 
$a_{i}(b)^{j-1}=c_{i}(b)^{j-1}=1$, put $c_{i}(b)^j=1$; for all other $i$ 
put $c_{i}(b)^j=0$.
\item Suppose first that the number of 1-bits in column $j$ is even. 
If 
\begin{equation}
\mbox{$a_{i}(b)^{j}\oplus c_{i}(b)^{j}=0$}
\end{equation} 
for every $i$, then pick $h$ and $l$ with $h\neq l$ such that
\begin{equation}
\mbox{$a_{h}(b)^{j-1}+c_{h}(b)^{j-1}=a_{l}(b)^{j-1}+c_{l}(b)^{j-1}=1$}, 
\end{equation}
and put $c_{h}(b)^{j}=c_{l}(b)^{j}=1$. $\{$We'll see later that such $h$ and 
$l$ indeed always exist.$\}$
\item Secondly, suppose that the number of 1-bits in column $j$ is odd. 
\begin{enumerate}
\item If for every $i$ for which $a_{i}(b)^{j-1}\oplus c_{i}(b)^{j-1}=0$ we
have $a_{i}(b)^{j}\oplus c_{i}(b)^{j}=0$, then pick $h$ such that
$a_{h}(b)^{j-1}+c_{h}(b)^{j-1}=1$ and 
$a_{h}(b)^{j}+c_{h}(b)^{j}=0$, and put $c_{h}(b)^{j}=1$.
\item If there is $i$ for which $a_{i}(b)^{j-1}\oplus c_{i}(b)^{j-1}=0$, and
$a_{i}(b)^{j}\oplus c_{i}(b)^{j}=1$, then pick $h$ such that
$a_{h}(b)^{j-1}+c_{h}(b)^{j-1}=1$ and put $c_{h}(b)^{j}=1$.
\end{enumerate}
\end{enumerate}
\end{enumerate}

{\bf Validity Proof of the Algorithm}\medskip

We begin by observing the general structure of the algorithm. In 
step 2(a) column $j$ of $c_1(b),\dots,c_n(b)$ is constructed. This
construction is consistent with Table~1. 
If we next go to step 2(b), then a correction to two of the 
carry bits might be done, by changing them from 0 to 1; if we 
go to step 2(c) instead, then a single carry bit will be changed 
from 0 to 1. No further corrections are done in column $j$. 

It suffices to show that the algorithm produces carry integers 
$c_1,\dots,c_n$ such that A, B, C of the above guidelines are 
satisfied. We will do this by showing that they hold for every 
column  $j$. This is clear for $j\leq k$ by step 1. In particular, 
C holds for $j=k$, since the multiset $S$ is balanced. (This is the only 
place in the proof where we use the fact that $S$ is balanced.) For 
$j\in\{k+1,\dots,m+1\}$ we use induction on $j$. So suppose A, B, C 
hold for column 
$j-1$ $(j\ge k+1)$, and we now apply the algorithm for column $j$. 

After applying step 2(a), which is consistent with the carry rules, 
suppose first that the number of 1-bits in 
column $j$ is even. We then say that column $j$ has {\it even parity}. 
If there is $h$ such that $a_{h}(b)^{j}+c_{h}(b)^{j}=1$, then there 
is also $l\ne h$ with $a_{l}(b)^{j}+c_{l}(b)^{j}=1$, since column $j$ 
has even parity, so property C holds. Otherwise, (5) holds for every 
$i$, and so C is violated. Now $h$ and $l\ne h$ with property (6) 
exist by the induction hypothesis. Moreover, in step 2(a) we have 
put $c_h(b)^j=c_l(b)^j=0$. Also $a_h(b)^j=a_l(b)^j=0$ by (5). So 
putting $c_h(b)^j=c_l(b)^j=1$ restores property C; it also preserves 
the even parity of column $j$, and is consistent with the carry rules.\medskip

We now suppose that, after applying step 2(a), column $j$ has 
{\it odd parity\/}, i.e., it contains an odd number of 1-bits, 
so step 2(c) applies. 

We assume first that the hypothesis of 2(c)i is satisfied. By the 
induction hypothesis, there is an {\sl even positive\/} number of 
$i$ for which $a_i(b)^{j-1}+c_i(b)^{j-1}=1$. For all of these $i$ 
we have $c_i(b)^j=0$ by step 2(a). Since column $j$ has odd parity, there 
thus exist $h$ and $l$ satisfying (6), for which, say, $a_h(b)^j+c_h(b)^j=0$ 
and $a_l(b)^j+c_l(b)^j=1$. Hence putting $c_h(b)^j=1$ restores both 
A and C, and is consistent with B. 

Secondly, assume that the hypothesis of 2(c)i is violated. Then the 
hypothesis of 2(c)ii holds. So there is $i$ for which 
$a_i(b)^{j-1}\oplus c_i(b)^{j-1}=0$, and 
\begin{equation}
\mbox{$a_i(b)^j+c_i(b)^j=1$}. 
\end{equation}
Note that $a_u(b)^{j-1}+c_u(b)^{j-1}=1$ implies $a_u(b)^j+c_u(b)^j\le 1$, 
since $c_u(b)^j=0$ by 2(a). 

\quad (I)\ Suppose that there is only a single $i$ for which 
$a_i(b)^{j-1}\oplus c_i(b)^{j-1}=0$ such that (7) holds. Since column 
$j$ has odd parity, the number of $u$ for which $a_u(b)^j+c_u(b)^j=1$ 
and $a_u(b)^{j-1}+c_u(b)^{j-1}=1$ must be even. Hence putting 
$c_u(b)^j=1$ for any such $u$ restores A and is consistent with C. 
Indeed any such $u$ is distinct from $i$, since 
$a_i(b)^{j-1}\oplus c_i(b)^{j-1}=0$, whereas $a_u(b)^{j-1}+c_u(b)^{j-1}=1$. 

\quad (II)\ Suppose there are at least two $i$ for which 
$a_i(b)^{j-1}\oplus c_i(b)^{j-1}=0$ such that (7) holds. Then C is 
already satisfied, so putting $c_u(b)^j=1$ for any $u$ as in (I) 
restores A and doesn't spoil C. 

Note that putting $c_u(b)^j=1$ in both (I) and (II) is consistent 
with B.\quad$\Box$\medskip


We now return to the proof of Theorem~1. It only remains to construct 
the $d_i$, which is done by the following algorithm. 
\begin{equation}
\mbox{For every $i=1$ to $n$ do: for $j=0$ to $m$ put $d_i(b)^j=c_i(b)^{j+1}$}.
\end{equation}
In other words, $d_i(b)$ is a ``right shift'' of $c_i(b)$. 

The validity of (8) is an immediate conclusion of Table~1 and Algorithm 
NotNimdi1: Table~1 shows that $d_i(b)^j=c_i(b)^{j+1}$ holds except for the 
second and penultimate rows. But when $a(b)^j=c(b)^j=0$, there is no 
reason to put $d(b)^j=1$, and when $a(b)^j=c(b)^j=1$, we may as well 
put $d(b)^j=1$. Thus these two rows do not arise in our case. (They 
may arise in the proof of Theorem~2, which follows below.)\quad$\Box$\medskip

{\bf Example:}
Let $\{a_{1},a_{2},a_{3}\}=\{3,5,8\}$. This is clearly a balanced multiset 
(with $k=0$). 

We have $3(b)=0011$ in the standard representation of binary numbers. 
Similarly, $5(b)=0101$ and $8(b)=1000$.

Following the steps of algorithm  NotNimdi1 we get
$c_{1}(b)=0100$, $c_{2}(b)=1010$, $c_{3}(b)=0000$. From (8), 
$d_{1}(b)=0010$, $d_{2}(b)=0101$, $d_{3}(b)=0000$, so 
$d_{1}=2$, $d_{2}=5$, $d_{3}=0$.

Since $d_{i}=t_{i}-a_{i}$ we have $t_{1}=5$, $t_{2}=10$, $t_{3}=8$. 
In binary, $t_{1}(b)=0101, t_{2}(b)=1010, t_{3}(b)=1000$.\par From 
all this we get that $t_{1}(b)\oplus t_{2}(b)\oplus t_{3}(b)=0111=7(b)$,
which is the same as
$(t_{1}(b)-a_{1}(b))\oplus (t_{2}(b)-a_{2}(b))\oplus (t_{3}(b)-a_{3}(b))$.
\bigskip

\textbf{Proof of Theorem~2.} We first show that the conditions are 
necessary. If $S$ is imbalanced, then even (1) doesn't hold, by 
Theorem~1. So suppose $S$ is balanced but not smooth, $n$ even, but 
$a_i(b)^k=1$ for all $i$. We have $c_i(b)^k=0$ for all $i$. Since $S$ 
is not smooth, there is an odd number of $a_i(b)^{k+1}=1$. To satisfy 
(4), we need an odd number of $c_i(b)^{k+1}=1$. This holds if and only 
if there is an odd number of $d_i(b)^k=1$, if and only if (2) is 
violated (since $d_i=t_i-a_i$). 

Finally, if $S$ is smooth but not even and $n=2$, then there is a 
least column $j$ such that $a_1(b)^j\oplus a_2(b)^j=0$ and 
$a_1(b)^{j+1}+a_2(b)^{j+1}=1$. 

We first consider the case where $a_1(b)^j=a_2(b)^j=1$. If also 
$c_1(b)^j=c_2(b)^j=1$, then $c_1(b)^{j+1}=c_2(b)^{j+1}=1$, so 
column $j+1$ has odd parity. The other possibility consistent 
with (4) is $c_1(b)^j=c_2(b)^j=0$. Then column $j+1$ has even 
parity if and only if $d_1(b)^j+d_2(b)^j=1$, and the latter 
contradicts (2). 

Secondly, let $a_1(b)^j=a_2(b)^j=0$. If $c_1(b)^j=c_2(b)^j=1$, 
then again column $j+1$ has even parity if and only if $d_1(b)^j+d_2(b)^j=1$. 
If $c_1(b)^j=c_2(b)^j=0$, then $c_1(b)^{j+1}=c_2(b)^{j+1}=0$, 
so column $j+1$ has odd parity.\medskip

For proving the sufficiency, we first consider the case where $n$ is 
odd. Since $S$ is balanced, Theorem~1 implies that there are integers 
$t_{1},\ldots,t_{n}$ with $t_{i}\geq a_{i}$, $i\in\{1,\ldots,n\}$, 
such that (1) holds. Let $d_{i}=t_{i}-a_{i}$. Then (3) holds. 

If $\sum_{i=1}^{'n} d_{i}(b)\not=0$, then there exists $j\in \{0,\ldots,m\}$
such that $\sum_{i=1}^{'n} d_i(b)^j=1$. This means that in the $n\times (m+1)$ 
matrix consisting of $d_1(b),\dots,d_n(b)$, the $j$-th 
column has odd parity. We wish to make it even, while, at the same time, 
preserving (3). 
 
At the beginning of the proof of Theorem~1 we saw that (3) holds if and 
only if (4) holds. Note that changing $d_i(b)^j$ may change $c_i(b)^{j+1}$. 
Table~1 shows, however, that for fixed $a_i(b)^j$, $c_i(b)^j$, a change 
in $d_i(b)^j$ does {\sl not\/} change $c_i(b)^{j+1}$ if and only if 
\begin{equation}
a_{i}(b)^{j}\oplus c_{i}(b)^{j}=0. 
\end{equation}

So it suffices to show that for every $j\in \{0,\ldots,m\}$, 
there is $i\in \{1,\ldots,n\}$ for which (9) holds (because this enables 
to regulate the parity of $d_i(b)$ so as to satisfy (2)). 
If this is not so, then $c_{i}(b)^{j}+a_{i}(b)^{j}=1$ for every $i$. 
Since $n$ is odd, we then have, 
\mbox{$\sum^{'n}_{i=1} (a_{i}(b)^{j}\oplus c_{i}(b)^{j}) = 1$}, 
which contradicts (4). 

Thus, for every $j$ for which $\sum_{i=1}^{'n} d_i(b)^j=1$ there is $i$ 
for which we can change $d_i(b)^j$ leaving (4), and hence (3), intact.\medskip 

Secondly, we examine the case where $n$ is even. 

{\sc Case I.} $S$ is balanced but not smooth, and $a_i(b)^k=0$ for 
some $i$. In Case~II below we indicate the changes the argument 
requires for the case where $S$ is smooth and $n\ge 4$. 

We construct the 
$c_i(b)$ by Algorithm NotNimdi2 below. It is a small modification of 
Algorithm NotNimdi1. The idea of the proof is to show that something 
like (9) holds for all the relevant columns $j$. To do this, we add 
another requirement to the guidelines A, B, C, namely:\medskip 

\qquad\qquad D.\ \ For every $j\in\{k,\dots, m+1\}$ there are 
$s,t\in\{1,\dots, n\}$, $s\ne t$, such that
$a_s(b)^j\oplus c_s(b)^j=a_t(b)^j\oplus c_t(b)^j=0$. 
Algorithm NotNimdi2 will implement the four guidelines. 

If $S$ is smooth, $n\ge 4$ and $a_i(b)^k=1$ for all $i$, we require 
D to hold only for $j\ge k+1$. 

\begin{center}{\bf Algorithm\quad NotNimdi2}
\end{center}

Steps 1, 2(a), 2(c)i are as in Algorithm NotNimdi1. Steps 2(b) 
and 2(c)ii are expanded:\medskip 

\qquad\ 2(b)\ Suppose first that the number of 1-bits in column $j$ 
is even. If either (5) holds for every $i$ or $a_i(b)^j+c_i(b)^j=1$ for 
every $i$, then pick $h$ and $l$ with $h\neq l$ such that (6) holds, 
and put $c_h(b)^j=c_l(b)^j=1$.\medskip

\qquad\ 2(c)ii\ If there is $i$ for which 
$a_{i}(b)^{j-1}\oplus c_{i}(b)^{j-1}=0$, and
$a_{i}(b)^{j}\oplus c_{i}(b)^{j}=1$, then pick $h$ such that
$a_{h}(b)^{j-1}+c_{h}(b)^{j-1}=a_{h}(b)^j+c_{h}(b)^j=1$ and put 
$c_{h}(b)^{j}=1$. If there is no such $h$, then pick $h$ such that 
$a_{h}(b)^{j-1}+c_{h}(b)^{j-1}=1$, and put $c_{h}(b)^{j}=1$.\bigskip 

{\bf Validity Proof of the Algorithm}\medskip

The validity proof is as that of Algorithm NotNimdi1, with the following 
additions. 

In column $k$, C is satisfied since $S$ is balanced and $c_i(b)^k=0$ 
for all $i$. Also D holds there, since there is $i$ for which 
$a_i(b)^k=0$ by hypothesis, and since $n$ is even. Incidentally, we 
see that $n\ge 4$. 

Suppose that C and D both hold for column $j-1$. We show that they 
hold also for column $j$\ ($j\ge k+1$). 

We consider first the case where column $j$ has even parity. If (5) 
holds for every $i$, then clearly both C and D are satisfied by putting 
$c_h(b)^j=c_l(b)^j=1$. So suppose that $a_i(b)^j+c_i(b)^j=1$ for 
every $i$.  By the induction hypothesis, there are integers $h,l$ 
satisfying (6). 
In step 2(a) we put $c_h(b)^h=c_l(b)^h=0$. So putting $a_h(b)^j=a_l(b)^j=1$ 
results in D being satisfied, and C is also satisfied since $n\ge 4$. 

Now consider the case where column $j$ has odd parity. By the 
induction hypothesis, there exist $h,l$, $h\ne l$ satisfying (6), 
and there exist $s,t$, $s\ne t$, satisfying 
$a_s(b)^{j-1}\oplus c_s(b)^{j-1}=a_t(b)^{j-1}\oplus c_t(b)^{j-1}=0$. 
In case~2(c)i, C has been restored for column $j$ without any change 
in rows $s$ and $t$, so D holds by the hypothesis of 2(c)i. 

In case~2(c)ii, if there is $h$ such that
$a_{h}(b)^{j-1}+c_{h}(b)^{j-1}=a_{h}(b)^j+c_{h}(b)^j=1$, then putting 
$c_h(b)^j=1$ makes $a_h(b)^j\oplus c_h(b)^j=0$. Since $n$ is even and 
A has been restored, there exists an index $i\ne h$ for which also 
$a_i(b)^j\oplus c_i(b)^j=0$. If, on the other hand, for every $h$ for which 
$a_{h}(b)^{j-1}+c_{h}(b)^{j-1}=1$ we have $a_{h}(b)^j+c_{h}(b)^j=0$, 
then putting $c_h(b)^j=1$ for one of these $j$ still leaves some 
$i$ for which $a_i(b)^j+c_i(b)^j=0$. Again, since $n$ is even, there 
are actually two distinct such $i$.\quad$\Box$\medskip


{\sc Case II.} $S$ is smooth and $n\ge 4$. If $a_i(b)^k=0$ for some 
$i$, then Case~I applies. We may thus assume $a_i(b)^k=1$ for 
all $i$. If either $a_i(b)^{k+1}=0$ for all $i$ or $a_i(b)^{k+1}=1$ 
for all $i$, put $c_u(b)^{k+1}=c_v(b)^{k+1}=1$ for some $u\ne v$. Then 
both C and D are satisfied for $j=k+1$. In any other case we have 
$a_u(b)^{k+1}=0$ and $a_v(b)^{k+1}=1$ for some $u,v\in\{1,\dots, n\}$. 
Since $S$ is smooth, there is actually an even number of $h$ satisfying 
$a_h(b)^{k+1}=1$, so at least 2. Since $n$ is even, there is an even 
number of $s$ such that $a_s(b)^{k+1}=0$, so at least 2. Putting 
$c_i(b)^{k+1}=0$ for all $i$, we see that both C and D are 
satisfied for $j=k+1$. 

Though D is not satisfied for $j=k$, it is clear that there is 
an even number of $d_i(b)^k=1$. In fact this holds precisely for 
the two values $u$, $v$ for which we put $c_u(b)^{k+1}=c_v(b)^{k+1}=1$ 
above. Now since C and D hold for $j=k+1$, they also hold for all 
$j\in\{k+1,\dots, m+1\}$ by the same induction proof used in 
Case~1.\quad$\Box$\medskip

In the previous example we found that 
$d_{1}(b)\oplus d_{2}(b)\oplus d_{3}(b)=0111=7(b)$. For $j\in\{0,1,2\}$, 
(9) holds only for $i=3$. This leads to the new value $d_3=7$, $t_3=15$, 
with $d_{1}(b)\oplus d_{2}(b)\oplus d_{3}(b)=0$.\bigskip

{\bf Proof of Corollary 1.} In the proof of Theorem~2 we observed 
that $n\ge 4$ also for the case where $S$ is balanced and $a_i(b)^k=0$ 
for some $i$. So for $n=2$, (2) holds if and only if $S$ is 
even.\quad$\Box$\bigskip

\begin{centering}
{\bf References}
\end{centering}\medskip

1. [BCG1982] E.R. Berlekamp, J.H. Conway and R.K. Guy, {\it Winning
Ways\/} (two volumes), Academic Press, London, 1982.

2. [Con1976] J.H. Conway, {\it On Numbers and Games}, Academic Press,
London, 1976.

3. [Fra1996] A.S.\ Fraenkel, Scenic trails ascending from sea-level Nim 
to alpine chess, in: {\it Games of No Chance\/}, Proc. MSRI Workshop 
on Combinatorial Games, July, 1994, Berkeley, CA (R. J. Nowakowski, ed.), 
MSRI Publ.\ Vol.~29, Cambridge University Press, Cambridge, pp. 13--42, 1996.

4. [FrL1991] A.S.\ Fraenkel and M.\ Lorberbom, Nimhoff games, 
{\it J.\ Combinatorial Theory\/} (Ser.~A) {\bf 58} (1991) 1--25.

5. [FrO1998] A.S.\ Fraenkel and M.\ Ozery, Adjoining to Wythoff's 
game its $P$-positions as moves, {\it Theoret.\ Comput.\ 
Sci. $($Math Games\/$)$} {\bf 205} (1998) 283--296.

6. [Guy1991] R.K.~Guy, editor, Combinatorial Games, {\it Proc. Symp. 
Appl. Math.} {\bf 43}, Amer. Math. Soc., Providence, RI, 1991.

7. [Now1996] R.J.~Nowakowski, editor, Games of No Chance, Mathematical 
Sciences Research Institute Publications, Vol.~29, Cambridge University 
Press, Cambridge, U.K., 1996.
\end{document}